\documentclass{article}

\usepackage{arxiv}

\usepackage[utf8]{inputenc} 
\usepackage[T1]{fontenc}    
\usepackage{hyperref}       
\usepackage{url}            
\usepackage{booktabs}       
\usepackage{amsfonts}       
\usepackage{nicefrac}       
\usepackage{microtype}      
\usepackage{lipsum}		
\usepackage{graphicx}
\usepackage{natbib}
\usepackage{doi}

\usepackage{amsmath}
\usepackage{graphicx}
\usepackage{amssymb}
\usepackage{booktabs}
\usepackage{color}
\usepackage{epsfig,amsthm}
\usepackage{hyperref} 
\usepackage{subcaption}

\title{Option Dynamic Hedging Using Reinforcement Learning}

\author{ Cong ZHENG\\
        \texttt{czhengae@connect.ust.hk}\\
 Department of Mathematics\\
 The Hong Kong University of Science and Technology\\
        \And Jiafa He \\
        \texttt{jhebu@connect.ust.hk}\\
        Department of Mathematics\\
 The Hong Kong University of Science and Technology\\
        \And  Can YANG\\
        \texttt{macyang@ust.hk}\\
 Department of Mathematics\\
 The Hong Kong University of Science and Technology\\
}

\renewcommand{\shorttitle}{\textit{arXiv} Template}

\hypersetup{
pdftitle={A template for the arxiv style},
pdfsubject={q-bio.NC, q-bio.QM},
pdfauthor={David S.~Hippocampus, Elias D.~Striatum},
pdfkeywords={First keyword, Second keyword, More},
}

\begin{document}
\maketitle

\begin{abstract}
	This work focuses on the dynamic hedging of financial derivatives, where a reinforcement learning algorithm is designed to minimize the variance of the delta hedging process. In contrast to previous research in this area, we apply uncertainty estimation technology to measure the uncertainty of the agent's decision, which can further reduce unnecessary wear and tear in the hedging process and control model overconfidence that may lead to significant losses. Numerical experiments show the superiority of our strategy in Monte Carlo simulations and SP 500 option data.
\end{abstract}


\section{Introduction}

\subsection{Background}
Hedging is a common practice used in trading to manage the risks associated with derivative transactions. Derivatives are financial instruments whose value is derived from an underlying asset, such as stocks, bonds, commodities, or currencies. Hedging involves taking an opposing position in a related asset or security to minimize or eliminate the risk of unfavorable price movements. One popular type of derivative is an option, which is a contract that allows the holder the right to buy or sell an asset at a specific price within a specific time period without any obligation.

In options trading, two parties are involved in a transaction: the option holder (or buyer) and the option writer (or seller). The option holder pays a premium to the option writer for the right to buy or sell the underlying asset at a predetermined price, known as the strike price. The option writer, in turn, assumes the obligation of the contract. This means that if the option holder chooses to exercise the option at expiration, the option writer must sell or buy the underlying asset to or from the option holder at the agreed-upon price. The option writer can also choose to buy back the option from the option holder before expiration to close out the position. Typically, option writers are institutional investors or market makers who earn option premiums by selling options, but they also bear the risk of price movements in the underlying asset. To mitigate this risk, option writers often use hedging strategies to limit their exposure to price volatility. For example, a market maker who sells call options on a stock could buy an equivalent amount of the underlying stock to hedge against the risk of a price increase.

Options dynamic hedging is the process of reducing risk by buying and selling the underlying asset. Its objective is to hedge the risk associated with selling options by buying or selling the underlying asset. This strategy is used by traders who wish to minimize the risk of their options position and maximize their profit potential. As the price of the underlying asset fluctuates, the risk exposure of the asset corresponding to the option also changes. To mitigate this risk, traders must maintain a dynamic hedge by continuously adjusting the quantity of the underlying asset they hold. For example, if a trader sells a call option on a stock and the stock price increases, the trader may face the risk of having to sell the stock at a lower price than the market value. To hedge against this risk, the trader may buy an equivalent amount of the stock to offset the potential loss. 

Reinforcement learning is a type of machine learning that focuses on learning the best action policy by having an agent interact with an environment to maximize long-term rewards. In the context of options dynamic hedging, reinforcement learning can be utilized to learn the best dynamic hedging strategy for maximizing profits in options trading. Reinforcement learning algorithms accomplish this objective by solving Markov Decision Processes (MDPs), where the environment is modeled as a sequence of states, actions, and rewards, and the agent makes decisions based on the current state. The agent receives feedback in the form of rewards from the environment, which it uses to optimize its decision-making strategy. By utilizing reinforcement learning to solve MDPs, traders can develop dynamic hedging strategies that adapt to changing market conditions and maximize their profits.

In options dynamic hedging, reinforcement learning can be utilized to learn the best dynamic hedging strategy. In this MDP problem, the state can include information such as the underlying asset price, exercise price, and current position holding size, while the actions are buying or selling the underlying asset. The reward is the profit generated from options and trading the underlying asset. The learning objective is to enable the agent to automatically adjust the position size of the underlying asset based on market changes and its own state, in order to minimize the risk exposure associated with selling options. By leveraging reinforcement learning, traders can develop dynamic hedging strategies that adapt to changing market conditions and maximize their profits.

\subsection{Related Work}
Traditional options hedging strategies rely on option pricing models, with the most widely used pricing model being the Black-Scholes model \cite{Black1973pricing}, also known as the BSM model. The optimal delta value, or the optimal position size of the underlying asset, can be derived by taking the partial derivative of the option price with respect to the underlying asset price. However, in real-life trading scenarios, replicating the payoff of derivatives through a continuous balancing of cash and stock quantities is difficult due to several factors. Firstly, the volatility of the underlying asset is not constant, which means that price fluctuations do not follow a perfect geometric Brownian motion. Additionally, continuous hedging is impossible due to transaction costs \cite{Leland1985option} and slippage. To address this challenge, a discrete-time hedging pattern, such as hedging on an hourly or daily basis, can be adopted. However, this solution can introduce significant errors in the hedging process. Hence, the optimal dynamic hedging strategy involves balancing hedging errors and underlying asset trading slippage. Increasing the hedging frequency can reduce hedging errors but also increase hedging costs, while decreasing the hedging frequency can reduce hedging costs but also increase hedging errors. By using reinforcement learning to develop dynamic hedging strategies, traders can optimize their hedging frequency and minimize their risk exposure while maximizing their profits.

Real-life scenarios often come with many restrictions, posing significant challenges to traditional theoretical models and making it difficult to obtain closed-form solutions through modeling. This complexity has a significant impact on the performance of dynamic hedging strategies. To address this challenge, \cite{Halperin2017qlbs} was the first to use reinforcement learning models to solve the problem of dynamic hedging of options. They used Q-learning reinforcement learning models with a limited action space and an approximate finite state to achieve dynamic hedging with limited actions. To better approximate real trading environments, \cite{dynamicreplication} considered the issue of hedging costs. Based on a given pricing model, the agent can automatically adjust the position size of the underlying asset by taking specific actions in a specific environment while simultaneously considering transaction costs and hedging risk in the reward function $r_t=\delta \Pi_t-\frac{\lambda}{2} \operatorname{Var}\left[\delta \Pi_t\right]$. Hedging costs $\delta \Pi_t$ arise from the profits generated from selling options and trading stocks, while hedging risk is approximated using the variance of the profits. By considering both hedging costs and hedging risk in the reward function, traders can develop dynamic hedging strategies that adapt to real-life trading scenarios and maximize their profits.

With the further development of reinforcement learning, \cite{Du2020deep} utilized more advanced reinforcement learning algorithms, such as DQN \cite{DQN} and PPO \cite{Schulman2017proximal}, to learn the action space for dynamic hedging of options. In their setup, the state variables typically include the time to expiration, stock price, current position of the underlying asset, and strike price, represented as $s_t=\left(t, S_t, N_t, K\right)$. Another common approach is to use the ratio of the stock price and option price as one of the state variables, known as moneyness, which enables simultaneous training of option data with different strike prices during the training process. This is particularly important when training with real data. 

Although a discrete action space brings greater convenience to the learning process, it inevitably introduces significant errors. To address this challenge, \cite{Cao2021deep} was the first to adopt a continuous action space, using the Deep Deterministic Policy Gradient (DDPG) reinforcement learning algorithm to enable the continuity of state and action spaces. In addition, unlike other authors who put hedging risk in the reward function, the author used two Q-functions to simultaneously estimate the expected loss and the expected square loss, making the selection of the objective function more flexible and the training more stable.

Undeniably, reinforcement learning has brought greater flexibility to dynamic option hedging by allowing different goals to be achieved through setting different states and reward functions. However, a key issue in applying deep learning to solve practical problems is how to model uncertainty \cite{kendall2017uncertainties}. The black-box nature of deep learning models is a major obstacle to their widespread use in some critical areas, as extreme instruction can cause incalculable economic losses. In the financial field, this problem is particularly prominent. Therefore, modeling uncertainty is an important component in the process of applying reinforcement learning models to solve option dynamic hedging. By accounting for uncertainty, traders can develop more robust and effective dynamic hedging strategies that can adapt to uncertain market conditions and minimize their risk exposure while maximizing their profits.

Generally, uncertainty is divided into two parts: aleatory uncertainty and epistemic uncertainty. Aleatory uncertainty arises from noise in the environment and is the error generated by observation. For aleatory uncertainty, we can usually model it, but it cannot be reduced. In the problem of option hedging, the fluctuation of the underlying asset or stock price is unpredictable or noise in the environment. For reinforcement learning problems, aleatory uncertainty usually has three sources: randomness in rewards, randomness in observations, and randomness in actions \cite{lockwood2022review}. Epistemic uncertainty is usually called model uncertainty, which arises from the limitations of the current model's learning ability. For this type of uncertainty, we can improve it by increasing the amount of training data and enhancing the model's learning ability.

The application of reinforcement learning models to solve dynamic hedging problems in options trading faces two main challenges:
\begin{itemize}
    \item We need to balance the conflict between the accuracy of option hedging and the cost of hedging. Pursuing high accuracy in hedging can lead to significant transaction cost wear and tear, while insufficient hedging can result in significant risk exposure.
    \item Deep models often provide overly confident answers, which can lead to significant real-world losses when solving problems in the financial domain. Uncertainty estimation can be used to address the issue of model overconfidence.
\end{itemize}

\subsection{Our Work And Contribution}
In this paper, we present a Bayesian reinforcement learning framework that allows us to learn a mapping from input data to aleatoric uncertainty and incorporate epistemic uncertainty in the decision-making process. Through experiments on both simulated and real-world data, we demonstrate that our framework achieves better results. Our main contributions can be summarized as follows:

First, this work is the first to analyze both aleatoric and epistemic uncertainty in dynamic option hedging, which has important practical significance.

Second, we made adjustments to the Q function in reinforcement learning according to the estimation of uncertainty, allowing for more accurate and robust dynamic hedging strategies.

Finally, by modeling uncertainty and obtaining attenuation from reward uncertainty, we can reduce the impact of noise on the model learning process and improve model performance.

The following chapters are organized as follows. Section 2 introduces reinforcement learning algorithms. Section 3 explains how to apply them to solve dynamic option hedging problems and perform uncertainty analysis in reinforcement learning. Section 4 presents experimental results on simulated data and S\&P 500 option data. Conclusions and analyses are presented in Section 5.

\section{Reinforcement Learning Algorithm}
Reinforcement learning (RL) is a reward-driven approach aimed at guiding continuous decision-making problems. Unlike supervised learning, which improves performance based on labeled training data, RL is more like an unsupervised learning approach where model performance is improved through the interaction with the environment instead of a labeled training dataset. This section introduces the basics of reinforcement learning, including state-action value functions and two types of strategy optimization methods: value-based optimization and policy-based optimization. The most famous method for value-based optimization is the deep Q-network \cite{DQN}, while policy-based optimization mainly introduces the deep deterministic policy gradient \cite{DDPG} algorithm. By understanding the fundamentals of reinforcement learning, traders can develop more effective dynamic hedging strategies that adapt to real-life trading scenarios. 

\subsection{RL Basics}
Reinforcement learning is often viewed as a Markov decision process, where at time t, the agent responds with action $A_t$ based on the observed current state $S_t$. The environment then provides feedback for this action in the form of a reward $R_t$, and the agent transitions to the next state $S_{t+1}$. This Markov decision process is defined as a quadruple $<S_t, A_t, R_t, S_{t+1}>$. The goal of reinforcement learning is to maximize the cumulative long-term return $G_t$, which is defined as:
\begin{equation}\label{Reward of rl}
G_{t=0: T}=R(\tau)=\sum_{t=0}^{\mathrm{T}} \gamma^t R_t,
\end{equation}
here, $\tau$ represents the trajectory, T represents the endpoint time node, and $\gamma$ is a constant used to balance long-term and short-term returns. The value of gamma is usually between 0 and 1, indicating the relative importance of immediate and future rewards. In order to express the expected long-term returns, we introduce the state value function and the action value function. The value function $V(s)$ represents the expected return of state $s$. For example, if there are two different states $S_1$ and $S_2$ in the next step, we can use $V(S_1)$ and $V(S_2)$ to evaluate their values. The policy typically selects the state with the higher value function as the next step state. If the agent's actions are based on a certain policy $\pi$, we can express the value function as:
\begin{equation}\label{Value function}
V^\pi(s)=\mathbb{E}_{\tau \sim \pi}\left[R(\tau) \mid S_0=s\right].
\end{equation}
To express the expected return based on both state and action, we define the action value function $Q^\pi(s, a)$ as follows:
\begin{equation}\label{Actio value function}
\begin{aligned} Q^\pi(s, a) & =\mathbb{E}_{\tau \sim \pi}\left[R(\tau) \mid S_0=s, A_0=a\right] \\ & =\mathbb{E}_{A_t \sim \pi\left(\cdot \mid S_t\right)}\left[\sum_{t=0}^{\infty} \gamma^t R\left(S_t, A_t\right) \mid S_0=s, A_0=a\right]. \end{aligned}
\end{equation}
If we have obtained an estimate $Q^\pi(s, a)$ for the action value function, then the optimal policy is to select the action $a_t$ that maximizes the value of the action value function for a given state $s_t$. In other words, the optimal action $a^*$ is: $a^*=\arg \max _a Q^\pi(s, a).$ We can observe the following relationship between the value function and the action value function:
\begin{equation}
    V_\pi(s)=\mathbb{E}_{a \sim \pi}\left[Q_\pi(s, a)\right].
\end{equation}

\subsection{Deep Q Learning}
In order to maximize the long-term reward function, Q-learning was introduced as an offline reinforcement learning method by \cite{Qlearning}. The algorithm updates the Q value by continuously trying different actions and receiving rewards from the environment. Specifically, when the agent is in a certain state $S_t$, by trying action $A_t$, we can obtain feedback from the environment, denoted as $R_{t+1}$ in equation \ref{time differential update}. Then, based on this information, we can update the previous estimate of the Q value using a time differential method. By iteratively updating the Q value, we can obtain the optimal estimate of the Q value and determine the optimal action policy, thus achieving maximum long-term return.
\begin{equation}\label{time differential update}
    Q\left(S_t, A_t\right) \leftarrow Q\left(S_t, A_t\right)+\alpha\left(R_{t+1}+\gamma \max _{a \in \mathcal{A}} Q\left(S_{t+1}, a\right)-Q\left(S_t, A_t\right)\right).
\end{equation}

However, when the environment becomes complex and the state becomes high-dimensional and continuous, it becomes difficult to estimate the Q value by just maintaining a Q table. Deep Q Learning proposes to use deep neural networks to approximate the Q value function, with the current state as the input to the neural network and the output being the Q value for each action. During the training process, the neural network parameters are optimized by minimizing the temporal difference (TD) error in equation \ref{TD Error}. Compared with the traditional Q-learning algorithm, Deep Q Learning can better handle high-dimensional and continuous state spaces.
\begin{equation}\label{TD Error}
    \left(R_{t+1}+\gamma \max _{a \in \mathcal{A}} Q\left(S_{t+1}, a\right)-Q(S_t, A_t )\right)^2.
\end{equation}

\subsection{Deep Deterministic Policy Gradient}
The Deep Deterministic Policy Gradient \cite{DDPG} algorithm can be viewed as a combination of the Deterministic Policy Gradient \cite{DPG} algorithm and deep neural networks, or as an extension of the Deep Q Learning \cite{DQN} algorithm in continuous spaces. Although we can discretize the action space in many problems, this process inevitably introduces errors. Unlike Q-learning, the DDPG algorithm simultaneously builds a policy function (Actor) and a Q-value function (Critic). The Q-value function part is consistent with DQN and is updated using the temporal difference algorithm, while the policy function part is updated using the policy gradient algorithm with the goal of maximizing the corresponding Q-value function.

We can represent the policy function as $\pi(S|\theta)$, which is a function of the state $S$ with the function parameter $\theta$. It is also a deterministic policy function, which means that the action in each round is directly obtained from the formula $A_t = \pi(S_t|\theta)$, rather than being sampled from a random policy.

The update iteration process for the action-value function $Q(s_t,a_t)$ in the DDPG algorithm is consistent with Deep Q Learning, which is updated through the Bellman equation. Additionally, to maintain stable learning, the DDPG algorithm employs a target network $Q^{'}$. The target network is an independent network that generates the target for Q-learning, further improving the stability of the neural network. Moreover, the target network adopts an exponential decay method to synchronize parameters with the main Q-network. This process ensures that the generation of the target value is not affected by the latest parameter updates, reducing the possibility of non-convergence. The loss function is adjusted to equation \ref{Target network TD Error}, which takes the difference between the predicted Q-value and the target Q-value generated by the target network.

\begin{equation}\label{Target network TD Error}
    \left(R_{t+1}+\gamma  Q^{'}\left(S_{t+1}, \pi^{\prime}\left(S_{t+1} \mid \theta^{\pi^{\prime}}\right)\right)-Q(S_t, A_t )\right)^2.
\end{equation}

\section{Applications In Option Dynamic Hedging}

\subsection{Background Theory}
The Black-Scholes model is a famous model used to solve option pricing and dynamic hedging. It assumes that the market is perfectly efficient and there are no transaction costs or taxes. Additionally, it assumes that the stock price follows a geometric Brownian motion, with the volatility being constant and the logarithmic changes in price following a normal distribution \ref{geometric Brownian motion}:
\begin{equation}\label{geometric Brownian motion}
    \delta S_t=\mu S_t \delta t+\sigma S_t \delta W_t,
\end{equation}

where $S_t$ represents the price of a stock at time $t$, $\mu$ represents the annualized average growth rate of the stock price, $\sigma$ represents the annualized volatility of the stock, and $W_t$ represents Brownian motion, which is commonly used to simulate stock price fluctuations. $dW_t$ is the incremental change in the Brownian motion, which follows a normal distribution $N(0,dt)$, where $dt$ represents the time interval. This formula describes the changes in the stock price between time $t$ and $t+dt$.

Under these self-financing assumptions, the option price $C_t$ can be easily derived:
\begin{equation}\label{call option price}
    C_t=S_t \Phi\left(d_1\right)-e^{-r \tau} K \Phi\left(d_2\right),
\end{equation}
where $S_t$ is the stock price at time t, r is the interest rate, K is the strike price, and $\tau$ represents the time to maturity $\tau=T-t$. $\Phi$ is the cumulative distribution function.
\begin{equation}\label{cumulative distribution function}
    \begin{split}
        d_1 &=\frac{\ln \left(S_t / K\right)+\left(r+\sigma^2 / 2\right) \tau}{\sigma \sqrt{\tau}}, \\
        \quad d_2&=d_1-\sigma \sqrt{\tau}, \\
        \Phi(x)&=\frac{1}{\sqrt{2 \pi}} \int_{-\infty}^x e^{-\frac{y^2}{2}} d y. \\
    \end{split}
\end{equation}

The option price can be seen as a function of five variables, $C_t = f(S_t,\tau,r, K,\sigma)$. For the dynamic hedging problem, we are concerned with the adjustment of the hedge position, which is the problem of determining the value of $\delta = \partial C / \partial S$. With the above formula \ref{call option price}, we can further calculate the analytical solution of $\delta$:
\begin{equation}\label{delta}
    \delta_t^{B S}=\frac{\partial C_t}{\partial S_t}= \Phi\left(d_1\right).
\end{equation}

It is important to be particularly attentive to this formula because the calculation of $\delta$ is the core of the problem of option dynamic hedging. Under the assumptions of the Black-Scholes model, we can analytically calculate the price of the option by modeling the changes in the stock price and then solve for the value of $\delta$. However, it is essential to note that in the process of obtaining this delta value, we have made several assumptions, including those regarding the stock price and continuous hedging. 

\subsection{RL Environment For Dynamic Hedging}

We consider the hedging problem of a call option. Suppose that at time t=0, in order to earn the premium, we sell a European call option. However, in order to reduce the loss caused by the risk of a one-sided upward market, we need to hold some stocks. The number of stocks we need to hold is the problem we need to consider, and the problem of dynamic hedging of options is to provide a corresponding solution to this problem.

To simplify the modeling process, we assume that the risk-free interest rate is 0. Under this assumption, we can focus on the key variables that drive the option dynamic hedging problem. Specifically, we assume that the stock price at time $t$ is denoted as $S_t$, and the option price at time $t$ is denoted as $C_t$. We also introduce a stock holding variable, denoted as $N_t$, which represents the number of stocks that we hold at time $t$. Additionally, we use $B_t$ to denote the amount of cash that we have at time $t$, and we consider the transaction fee generated by trading stocks, which is denoted as $f(S_t,\delta N_t)$. Here, $\delta N_t=N_{t+1}-N_t$ represents the change in stock holding at time $t$, or the number of stocks traded at time $t$.

In most cases, the transaction fee function $f$ is a linear function of $S_t*\delta N_t$, which means that the transaction fee is positively correlated with the trading volume. However, some scholars \cite{dynamicreplication} have proposed the use of nonlinear transaction fees to account for market impact and trading wear.

We can evaluate the value of our entire investment portfolio at time $t$ by using the variable $\Pi_t$. This portfolio value can be broken down into three components: cash, stocks, and options.
\begin{equation}\label{portfolio value}
    \Pi_t = B_t+S_t*N_t-C_t.
\end{equation}

If we achieve perfect hedging of the option, it means that the value of our investment portfolio does not change over time. This implies that the change in the value of the stock is used to offset the risk of selling the option. To achieve this perfect hedging, traders typically adjust the number of stocks they hold at time $t$ to hedge against the potential risk of the option. Suppose that the number of stocks changes from $N_t$ to $N_{t+1}$ at time $t+1$:
\begin{equation}\label{portfolio value change}
    B_{t+1} = B_t-S_t*(N_{t+1}-N_t)-f(S_t,\delta N_t) .
\end{equation}
By comparing the value of our investment portfolio at times $t$ and $t+1$, we can determine the change in the value of the portfolio over time. This change reflects the cost of dynamic hedging, which is the cost of adjusting the number of stocks we hold to hedge against the potential risk of the option. By tracking the magnitude of this cost, we can better understand the effectiveness of our hedging strategy:
\begin{equation}\label{Reward}
    R_t = C_t-C_{t+1}+N_{t+1}(S_{t+1}-S_t)-f(S_t,\delta N_t) ,
\end{equation}
the variable $R_t$ (also referred to as the Reward, see Equation \ref{Reward}) represents the change in the value of our investment portfolio from time $t$ to $t+1$. Since this change reflects the profit of our dynamic hedging strategy, it is often used as the main reward function in reinforcement learning algorithms.

In addition to considering hedging costs, controlling risk is another critical aspect of developing a successful trading strategy. One commonly used approach to this problem is mean-variance optimization, which aims to maximize the hedging return $E(R_t)$ while minimizing the hedging risk $Var(R_t)$.
\begin{equation}\label{mean variance optimization}
\begin{aligned}
    \begin{split}
        \text{Maximize } & E(R_t) \
        \text{subject to } & Var(R) < \delta \
    \end{split}
\end{aligned}
\end{equation}
To find the optimal balance between hedging loss and hedging risk, traders can adjust the value of a parameter $\gamma$. This parameter controls the relative weight of the reward function, which is a critical element of the reinforcement learning model. We propose an integrated reward function that reflects both the hedging loss and hedging risk, as follows:
\begin{equation}\label{mean variance reward}
    R_t = E(R_t) + \gamma Var(R_t).
\end{equation}

In addition to designing a loss function, selecting appropriate states is another crucial consideration when developing reinforcement learning models. In the context of selling a European call option, traders must buy a certain amount of stocks to hedge against the risk of an upward market. By analyzing the cumulative distribution function (see Equation \ref{cumulative distribution function}), we can see that the optimal position size $\delta$ under the Black-Scholes model assumption is related to several key variables, including the current stock price $S_t$, option exercise price $K$, time to expiration $\tau$, and stock volatility $\sigma$.

While the Black-Scholes model provides a useful framework for option pricing, it has limitations when applied to real-world data. One significant issue is that the model assumes a constant implied volatility, while the volatility surface \cite{gatheral2011volatility} in real data suggests that this assumption is not always valid. To address this issue, we propose including an estimate of the implied volatility as a state variable in our reinforcement learning model. Specifically, we use the historical volatility of the asset as an estimate of the implied volatility, denoted as $\sigma^t$, where $t$ represents the volatility estimate at different time dimensions.

In addition to historical volatility, we also include other variables from the Black-Scholes model as state variables, including $\gamma$ and $\theta$. Here, $\gamma$ represents the sensitivity of option prices to changes in the underlying asset prices, while $\theta$ represents the first-order derivative of option prices with respect to time. By incorporating these variables into the state space of our reinforcement learning algorithm, we can develop more accurate and effective trading strategies that account for the complexities of real-world market data.

In other words, we divide the state variables into two main categories:

\begin{itemize}
    \item The first category includes the agent's own relevant variables, such as the current stock price, option exercise price, time to expiration, and current position size.
    \item The second category includes predictive variables for the stock, such as volatility variables at different lookback times and some option price-related predictive factors. These variables enable the agent to anticipate future market trends and adjust their trading strategies accordingly.
\end{itemize}

In this part, we introduce the Bayesian reinforcement learning framework we propose for the dynamic hedging problem of options. Our model provides uncertainty estimates while solving the dynamic hedging problem, which is critical for accurately assessing the risks. By modeling the aleatoric uncertainty (in this case, uncertainty in the reward function caused by the uncertainty in stock prices), we can reduce the impact of noise on the training process and make the training more robust.

We have observed that the timing of our hedging strategy results in a high correlation between the hedging error, as measured by the reinforcement learning reward function $R_t$, and the change in stock price $S_{t+1}-S_t$. As predicting changes in stock price is a challenging task, this poses a significant challenge to our training process. To address this, we make the assumption that the reward function follows a normal distribution:
\begin{equation}\label{reward distribution}
    R_t \sim N(E(R_t),\sigma_t^{2}) .
\end{equation}

The uncertainty $\sigma_{t}^2$ is derived from the volatility of stock prices, which can be highly unpredictable and can have a significant impact on the effectiveness of our trading strategy. To account for this, we have included the historical stock volatility as a factor in our state variables. However, this also means that we need to provide an additional estimate in our neural network to accurately model this factor and its impact on our trading decisions.

Using the DDPG network as an example, we can define the actor network as $\mu(s|\theta^u)$ and the critic network as $Q(s,a|\theta^Q)$. The actor network parameters are denoted as $\theta^u$, while the critic network parameters are denoted as $\theta^Q$. We also have corresponding target networks, denoted as $Q'$ and $\mu'$. To estimate the epistemic uncertainty, we modify the output of the actor network to include two parts: $a_t$ and $\sigma_{t}^{2}$. In other words, instead of just outputting the action $a_t$ for a given state $s$, the actor network will also output an estimate of the uncertainty (i.e., variance) associated with that action.

To model the aleatoric uncertainty arising from the normal distribution of the reward function, we use the log-likelihood of a Gaussian distribution. The loss function can be transformed as follows:
\begin{equation}\label{aleatoric unvcertainty TD error}
    L_{BNN}(\theta^Q) = \frac{1}{D} \sum_t \frac{1}{2} \hat{\sigma}_t^{-2}\left\|\mathbf{y}_t-\hat{\mathbf{y}}_t\right\|^2+\frac{1}{2} \log \hat{\sigma}_t^2,
\end{equation}

here, $\sigma_t^{2}$ represents the estimated uncertainty in the reward function. The modified loss function is derived from the log-likelihood of a Gaussian normal distribution. It consists of two key components: the first component ensures that the TD error tends towards zero, but reduces the weight of samples with high uncertainty (which can be interpreted as samples with high volatility), while the second component serves as a regularized term to prevent the uncertainty estimate from approaching infinity.

This type of model has been widely studied in the field of deep learning \cite{nealbayesian,dikovbayesian}. Bayesian neural networks no longer use deterministic network weights, but instead use a distribution to replace deterministic parameters. In deep learning models, assuming that the network parameters are $W \sim \mathcal{N}(0,I)$, by sampling the values of $W$, we can obtain different network outputs $f^W(x)$. From this, we can obtain the probability distribution of the output $p(y|f^W(x))$.

Dropout variational inference is a commonly used technique in the field of deep learning. Dropout is typically a regularization method that randomly removes some neurons or sets some weights to zero during training to reduce overfitting, but in general, no neurons are removed during testing. MC dropout was first proposed as an inference method by Gal and Ghahramani \cite{galbayesian}, used to estimate the uncertainty of the model. For a trained model, different network outputs can still be obtained by performing multiple forward passes during testing, which can be used to obtain the probability distribution of the output and calculate the uncertainty of the output. Therefore, MC dropout can be seen as an extension of the Dropout technique, used to improve the confidence and robustness of the model. 

For regression tasks, uncertainty can be simply represented by the variance of the output results.
\begin{equation}\label{Epistemic uncertainy}
    Var(y) = \frac{1}{n-1} \sum_{t=1}^{T} \mathbf{y_t}^{T}\mathbf{y_t} - \mathbf{\Bar{y}}^T\mathbf{\Bar{y}},
\end{equation}
where $\mathbf{y_t} = f^{\mathbf{W}_t}(\mathbf{x}), \mathbf{\Bar{y}} = \frac{1}{T} \sum_{t=1}^T \mathbf{f}^{\mathbf{W}_t}(\mathbf{x})$.

In reinforcement learning models, the Q value is typically considered as an estimate of the model's expected future returns. To reflect the model's estimation of the Q value, we use the MC dropout method. By using dropout, we can obtain multiple Q values, which can be used to estimate the model's uncertainty. Specifically, we can obtain an estimate of the Q value uncertainty as follows:
\begin{equation}\label{Epistemic uncertainy in rl}
    Var(Q) = \frac{1}{n-1} \sum_{t=1}^{T} \mathbf{Q_t}^{T}\mathbf{Q_t} - \mathbf{\Bar{Q}}^T\mathbf{\Bar{Q}},
\end{equation}
where $Q_t = Q(s,a|\theta^Q_t)$ and $\mathbf{\Bar{Q}} = \frac{1}{T} \sum_{t=1}^T Q(s,a|\theta^Q_t)$.

\begin{figure}
	\begin{center}
	\includegraphics[width=0.98\columnwidth]{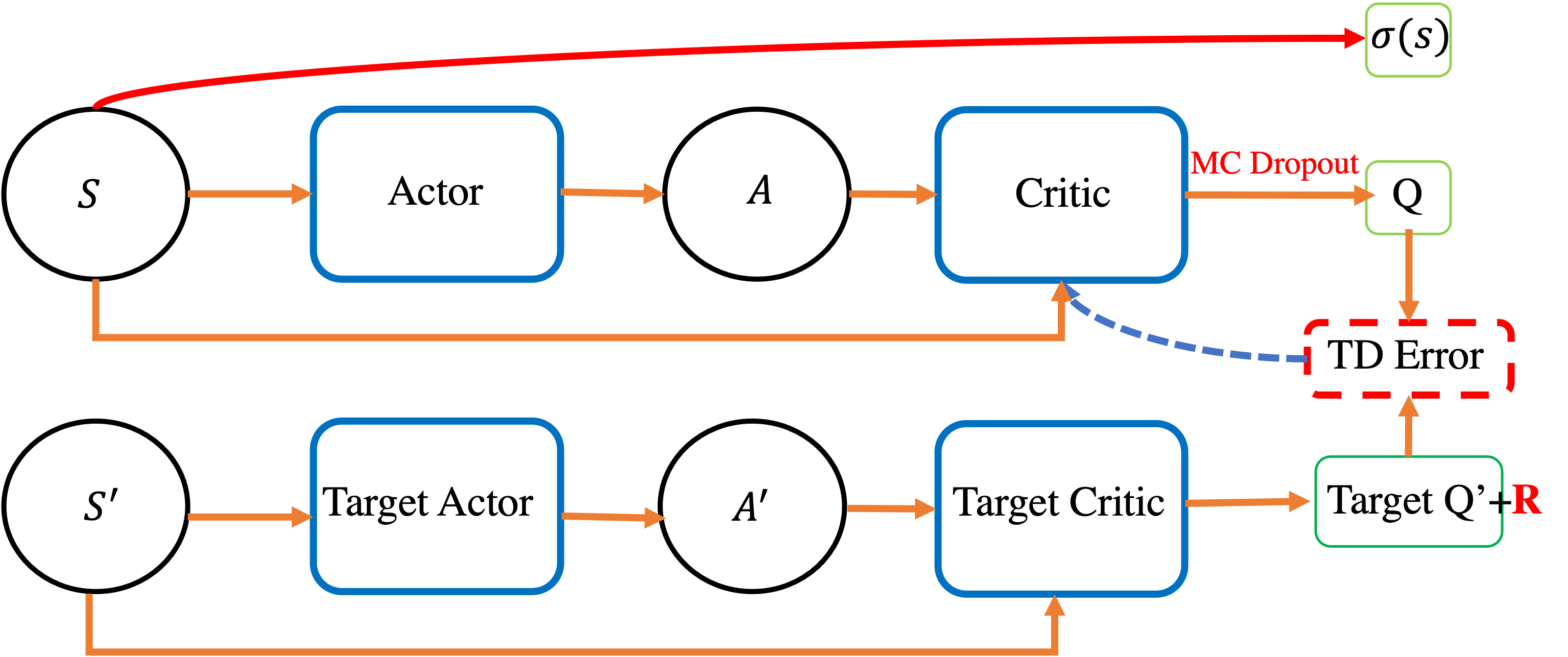
 }
	\end{center}
	\caption{\textbf{Uncertainty-aware DDPG network.} Our model architecture inherits from the Deep Deterministic Policy Gradient (DDPG) network. Different from the original network, we have added an output of $\sigma(t)$ to estimate the aleatoric uncertainty, and the critic network structure has been trained using the MC dropout method to estimate the epistemic uncertainty.}
	\label{fig1}
\end{figure}

Our model architecture inherits from the Deep Deterministic Policy Gradient (DDPG) network in \ref{fig1}.

The trading model is composed of two main parts: a local network and a target network. The local network consists of two sub-networks: the actor network and the critic network. The actor network is responsible for generating actions, while the critic network evaluates Q-values. The target network, depicted in the lower part of the figure, is used to calculate target values.

During training, the target network's parameters are primarily updated from the local network via soft updates. It's important to note that, unlike the original DDPG network, we have added an output of $\sigma(t)$ to estimate the aleatoric uncertainty, and the critic network structure has been trained using the MC dropout method to estimate the epistemic uncertainty.

Combining both aleatoric uncertainty and model uncertainty, the problem becomes how to minimize the equation:
\begin{equation}\label{Target network uncertainty aware TD Error}
    \frac{1}{2}\sigma(S_t)^{-2}\left(R_{t+1}+\gamma  Q^{'}\left(S_{t+1}, \pi^{\prime}\left(S_{t+1} \mid \theta^{\pi^{\prime}}\right)\right)-Q(S_t, A_t )\right)^2+\frac{1}{2}\log \sigma(S_t)^2.
\end{equation}

\section{Results}
The experimental section of the research paper is divided into two main parts: the simulation part and the real data solution part. In the simulation part, we provide some simple experimental results for option dynamic hedging problems. This includes generating Black-Scholes option pricing curves, theoretical delta values, and reinforcement learning solutions that take into account hedging costs. We also analyze the advantages of uncertainty in simulated scenarios. The real data solution part of the experiment demonstrates the differences between real data and simulated data, and how the proposed reinforcement learning model can be applied to solve the dynamic hedging problem of S\&P 500 options.

\subsection{Simulation Study}

\subsubsection{Experiment Results of Delta hedging}
In the simulation section, we assume that the stock data follows a geometric Brownian motion:
\begin{equation}\label{brownian motion}
    \delta S_t=\mu S_t \delta t+\sigma S_t \delta W_t.
\end{equation}

The experimental setup assumes that the parameter $\mu = 0.05$ represents the average annual growth of the stock, while the parameter $\sigma = 20\%$ represents the constant annual volatility of the stock at 20\%. These parameters are used to simulate the prices of the underlying stock. For the option data, the experiment assumes that we are dealing with at-the-money options that expire in one month, and the option prices follow the Black-Scholes pricing model. To hedge the risk associated with the options, the experiment considers timing hedging and artificially limits the hedging costs. This means that the model is designed to make trading decisions based on specific time intervals, rather than continuously adjusting the hedging position. For transaction fees, we follow the setup in \cite{Cao2021deep} and consider a $1\%$ transaction fee.

\begin{figure}
	\begin{center}
	\includegraphics[width=0.98\columnwidth]{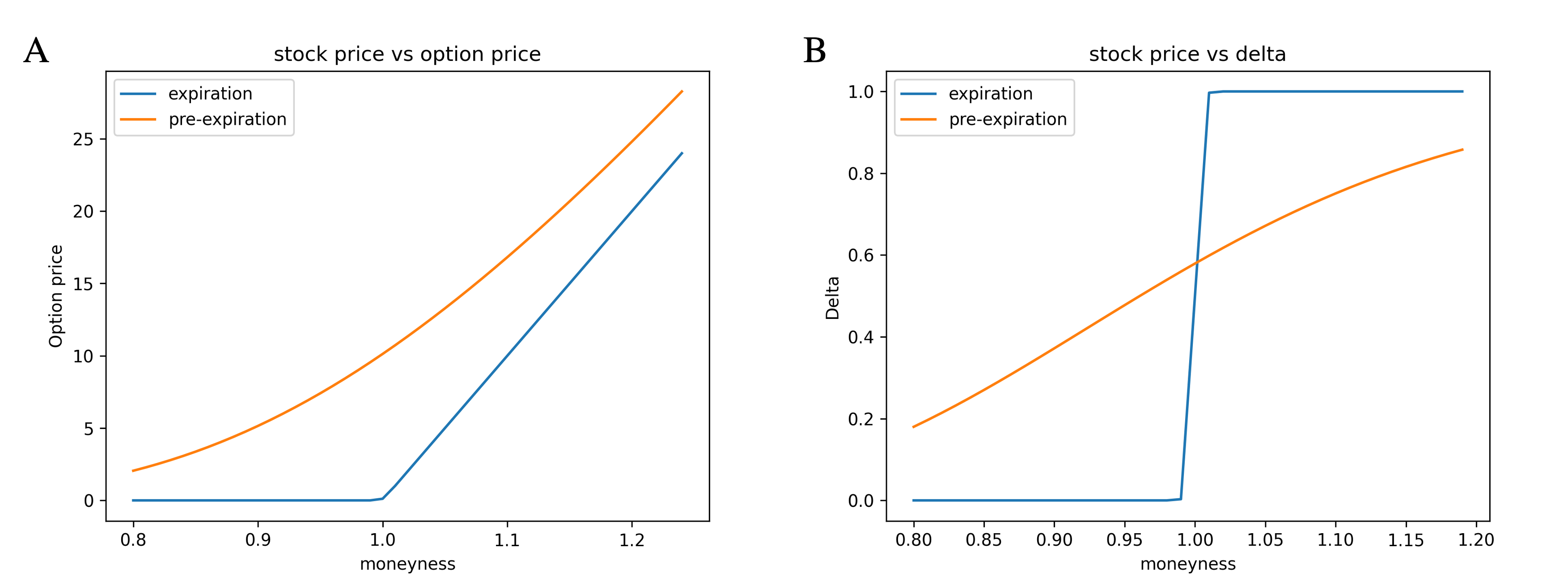}
	\end{center}
	\caption{\textbf{Option price and theoretic delta curve.} \textbf{A}: Relationship between stock price and option price: The blue curve represents a scenario in which the time to expiration is close to zero, while the yellow curve represents a scenario in which the time to expiration is far away. \textbf{B}: Relationship between moneyness and theoretical delta}
	\label{fig2}
\end{figure}

We first analyze the relationship between option prices, time to expiration, and stock prices. Assuming an initial stock price of $S_0 = 100$ and an option strike price of $K=100$, we define the variable "moneyness" as $moneyness = S_t/K$. For an at-the-money option expiring in one month, Figure \ref{fig2}A reflects the relationship between stock price and option price. In the graph, there are two lines: the blue line represents the option price when the time to expiration is close to zero, and the yellow line represents the option price when the time to expiration is far away. When the time to expiration is close to zero, the option price is composed of two parts: when the stock price is lower than the strike price, the option price tends to zero, and when the stock price is higher than the strike price, the option price equals its intrinsic value. When the time to expiration is far away, the option value equals its intrinsic value plus its time value. Another observation from the graph is that when the time to expiration is the same, the higher the stock price, the higher the option price. This is because a higher stock price increases the likelihood that the option will be exercised, which in turn increases the value of the option.

\begin{figure}
	\begin{center}
	\includegraphics[width=0.8\columnwidth]{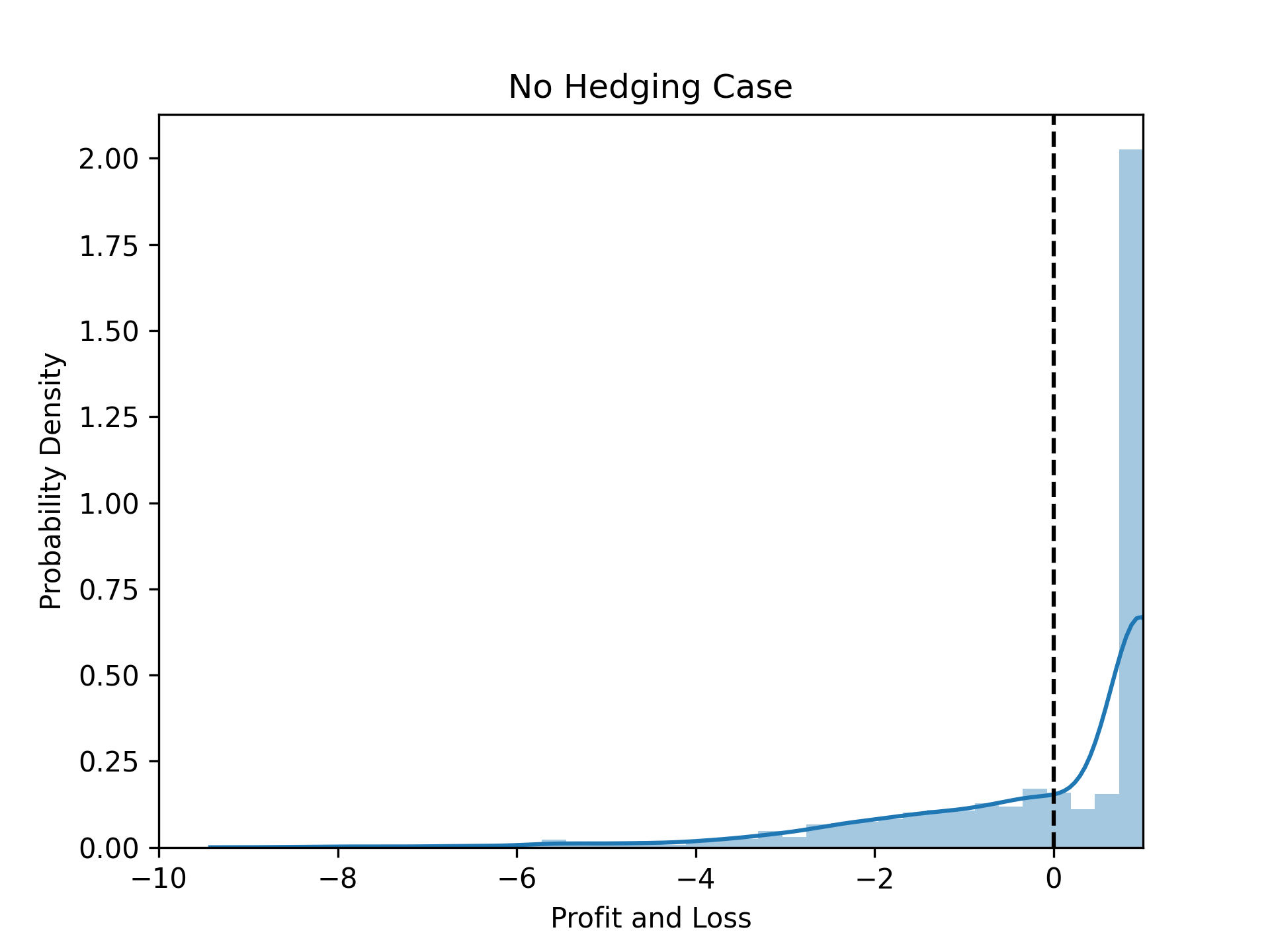}
	\end{center}
	\caption{\textbf{Distribution of returns of selling call option.} This figure reflects the expected profit graph of the call option seller with no hedge case, where E(R) = -4\%, Var(R) = 156\%. }
	\label{fig3}
\end{figure}

The second graph in the experimental analysis illustrates the relationship between stock prices and theoretical delta values. To hedge the risk exposure from selling options, it's common to buy $\delta$ shares of stock to hedge the risk. Based on the relationship curve between option prices and stock prices, we can obtain the $\delta$ value that satisfies the equation $dC = \delta dS$, which represents the number of shares of stock that we need to purchase and is also the slope of the stock option price curve. Graph B in Figure \ref{fig2} shows two curves. The blue curve represents the case where the time to expiration is 0. We can see that when the stock price is greater than the strike price, the delta value tends to 1, and when the stock price is less than the strike price, the delta value tends to 0. The yellow delta curve reflects the case where the time to expiration is not equal to 0 and is a gradually changing curve from 0 to 1 as the stock price increases.

\begin{figure}
	\begin{center}
	\includegraphics[width=0.8\columnwidth]{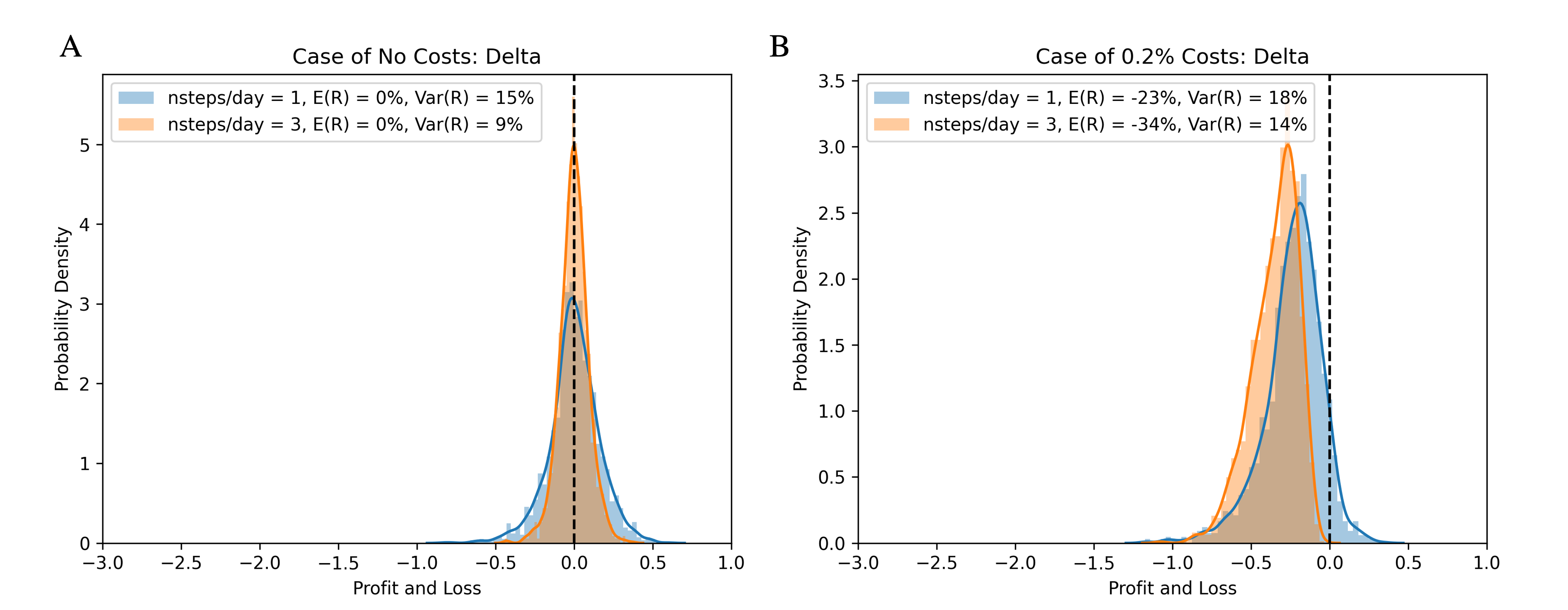}
	\end{center}
	\caption{\textbf{Distribution of returns of delta hedging.}  \textbf{A}: Hedging without considering the cost of buying and selling stocks. \textbf{B}: Hedging considering the cost of buying and selling stocks.}
	\label{Hedging with and without cost}
\end{figure}

We first consider the case of a seller of options who does not employ dynamic hedging. Since the price follows a Brownian motion, they calculate the expected return and variance of the seller's profit, represented by $E(R)$ and $Var(R)$, respectively.For example, when the initial stock price is 100, the time to expiration is 30 days, the annualized volatility is 20\%, and the strike price is also 100, the option price is approximately 2.28. To reflect the hedging effect, the authors consider the hedging cost as a percentage of the option price.

Figure \ref{fig3} reflects the expected profit graph of the call option seller. The graph shows that most of the profits are concentrated at the 100\% position, which is because most options have not reached the strike price, and the initial selling price of the option is the profit of the option seller. However, we also note that in a small number of scenarios, selling options can lead to catastrophic losses, with losses of up to 1000\%. This usually occurs in extreme market conditions where options are sold but the stock price keeps rising. In other words, the variance of returns is very high, and the risk is high, which is extremely detrimental to the strategy. This is also one of the reasons why dynamic hedging is needed.

To reduce the risk associated with selling options, we can buy a corresponding amount of stocks to implement dynamic hedging, so as to gain profits from the stocks when the stock price rises and offset the losses caused by selling options. First, we consider the dynamic hedging strategy $\delta$ given by the theoretical model, Black-Scholes model.

The simulation results \ref{Hedging with and without cost} are divided into two parts: hedging without considering the cost of buying and selling stocks and hedging considering the cost of buying and selling stocks. First, let's look at the left graph, which contains two sets of comparative experiments representing the results of different hedging frequencies. The blue color represents the distribution of daily hedging returns, while the yellow one represents the distribution of returns with hedging three times a day. It can be seen that the expected returns of both tend to zero. Furthermore, as the hedging frequency increases, the variance of returns decreases. Compared with the variance of $156\%$ without hedging, dynamic hedging significantly reduces the risk of the strategy, with variances of $9\%$ and $15\%$, respectively. The right graph shows the results considering the transaction cost. The simulation results show that as the hedging frequency increases, the expected return decreases from $-23\%$ to $-34\%$, while the variance is controlled to a certain extent, decreasing from $18\%$ to $14\%$.

\begin{figure}
	\begin{center}
	\includegraphics[width=0.98\columnwidth]{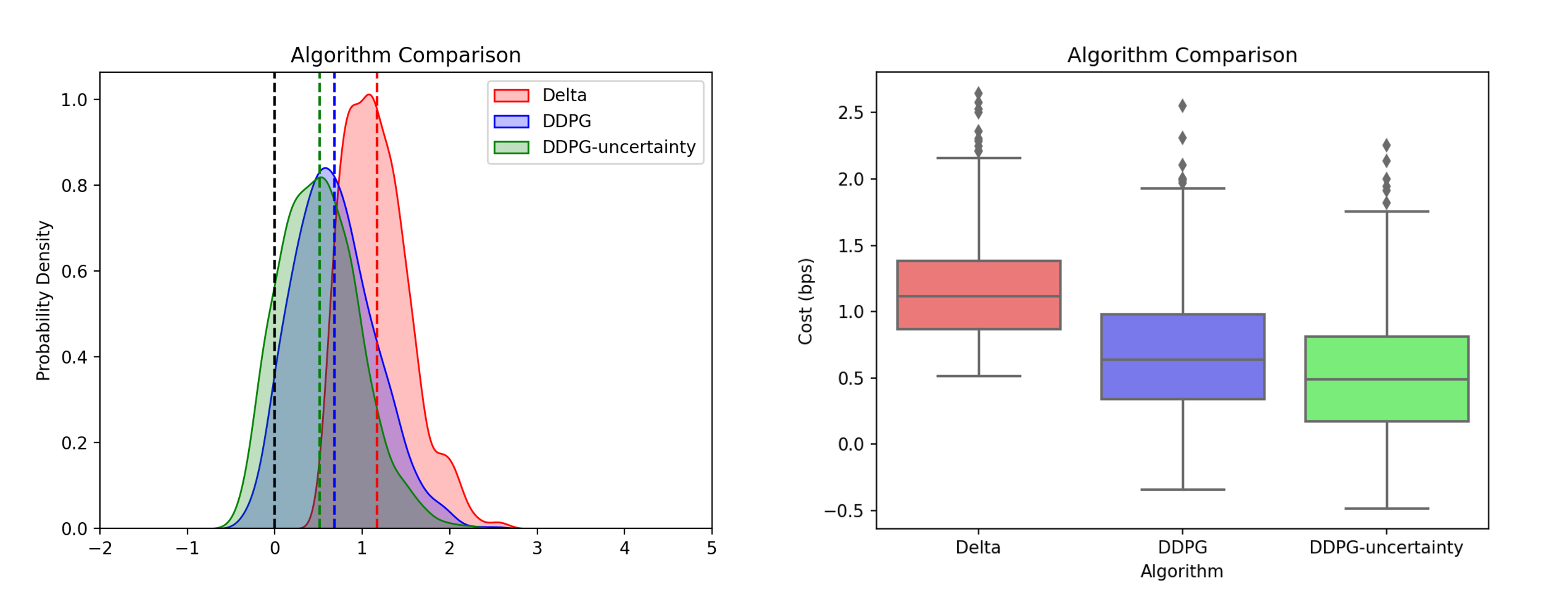}
	\end{center}
	\caption{\textbf{Strategy comparison chart.} This figure displays the probability density function of the strategy returns, comparing the theoretical value $\delta$ of option dynamic hedging with two models trained through reinforcement learning: DDPG without considering return uncertainty, and DDPG uncertainty considering return uncertainty.}
	\label{Strategy comparison}
\end{figure}

The simulation results show that when there is no hedging cost, the expected return is approximately equal to zero, and the higher the hedging frequency, the smaller the variance of the return. However, after considering factors such as transaction costs, the expected value of hedging becomes negative and the variance increases. The $\delta$ value provided by the Black Scholes model is no longer the optimal choice.

\subsubsection{Experiment Results of RL model}
Data-driven methods provide a more flexible solution for option dynamic hedging in complex real-world scenarios. In the experimental analysis, we compare the original reinforcement learning algorithm \ref{Target network TD Error} with our improved training method considering the uncertainty of returns \ref{Target network uncertainty aware TD Error}. By considering the expected returns and variance of the dynamic hedging returns of the two algorithms, we hope to evaluate their effectiveness in achieving higher expected returns and lower variance. 

\begin{table}[]
    \centering
    \begin{tabular}{|c|c|c|c|}
    \hline
     & Delta & DDPG & DDPG-uncertainty \\
    \hline
    mean & 118\% & 73\% & 58\%\\
    \hline
    std & 39\% & 42\% & 42\%\\
    \hline
    \end{tabular}
    \caption{\textbf{Strategy comparison table.} This table provides specific performance indicators, which compare the performance of the three strategies in terms of expected returns and variance.}
    \label{Strategy comparison table}
\end{table}

Figure \ref{Strategy comparison} in the experimental analysis displays the probability density function of the strategy returns, comparing the theoretical value $\delta$ of option dynamic hedging with two models trained through reinforcement learning: DDPG without considering return uncertainty, and DDPG uncertainty considering return uncertainty. The results show that the reinforcement learning models generally achieve better expected returns than the theoretical value of $\delta$, while controlling the risk of the strategy. Considering the uncertainty further improves the expected returns of the strategy. We also provide specific performance indicators in Table \ref{Strategy comparison table}, which compare the performance of the three strategies in terms of expected returns and variance.

\subsubsection{Results analysis}

\begin{figure}
	\begin{center}
	\includegraphics[width=0.8\columnwidth]{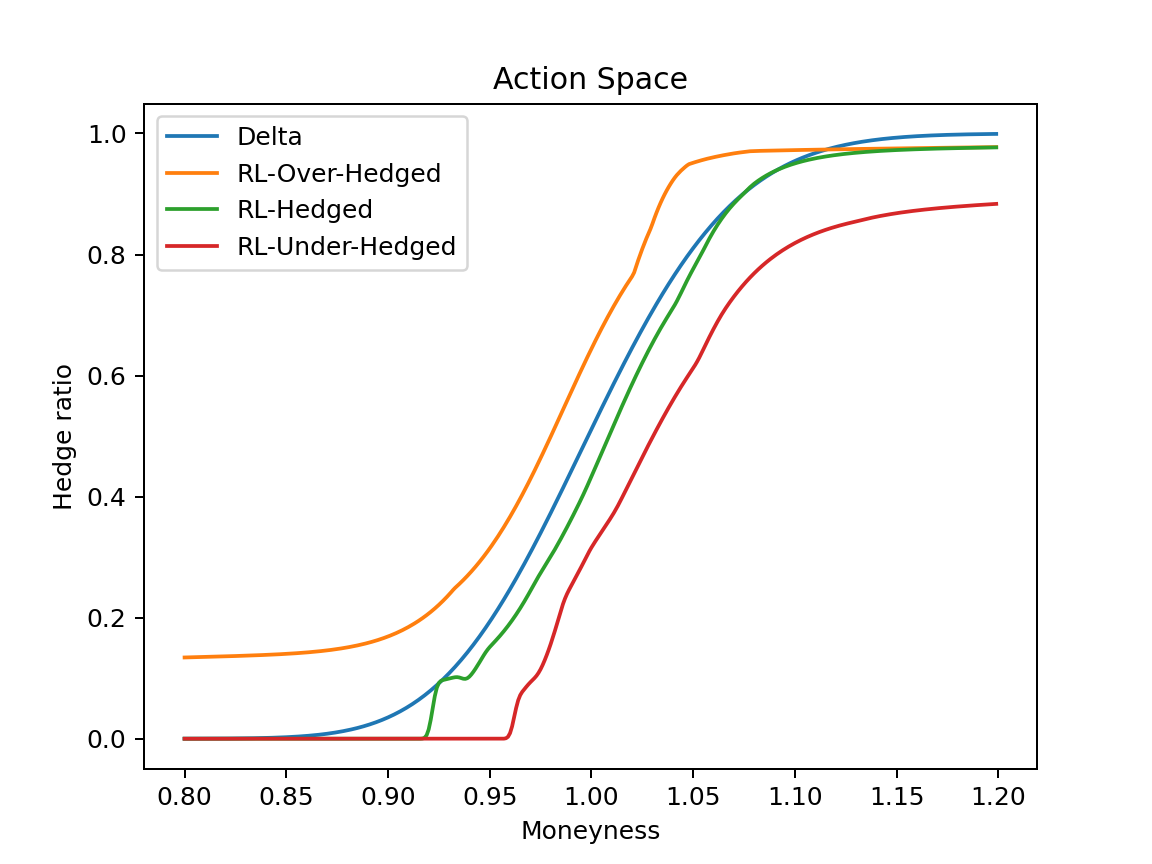}
	\end{center}
	\caption{\textbf{Agent action pattern.} This figure illustrates the rough decision logic of the agent in a dynamic hedging strategy. The horizontal axis represents the moneyness, which is the ratio of the stock price to the strike price, and the vertical axis represents the value of the variable $\delta$.}
	\label{Agent action pattern}
\end{figure}
The experimental results in the analysis show that although delta hedging provides the optimal solution in the no-cost case, in the presence of transaction costs, the reinforcement learning model provides a better solution. When implementing dynamic hedging strategies, it's important to consider the action space of the agent. The action space refers to the set of possible actions that the agent can take in response to market conditions.

Figure \ref{Agent action pattern} in the experimental analysis illustrates the rough decision logic of the agent in a dynamic hedging strategy. The horizontal axis represents the moneyness, which is the ratio of the stock price to the strike price, and the vertical axis represents the value of the variable $\delta$. The variable $\delta$ takes values from 0 to 1, with higher values indicating a higher proportion of the underlying asset to be held to replicate the option's payoff. As the stock price rises, the probability of exercising the option increases, causing the value of delta to approach 1. 

The figure compares the difference between the theoretical delta value and the delta value given by reinforcement learning. It can be seen that if the current stock position is approximately equal to the theoretical delta value, the decisions given by the two (the blue and green curves in the figure) are relatively close. However, if the current stock position is greater than the theoretical delta value (the yellow curve in the figure), the optimal decision at this point is not the theoretical delta value, but should be slightly larger than the theoretical value. Similarly, if the current stock position is less than the theoretical delta value (the red curve in the figure), the optimal decision at this point is not the theoretical delta value, but should be slightly smaller than the theoretical value. The decision given by reinforcement learning fully considers the current position and reduces unnecessary trading by adjusting the hedging position based on market conditions.

\begin{figure}
	\begin{center}
	\includegraphics[width=0.95\columnwidth]{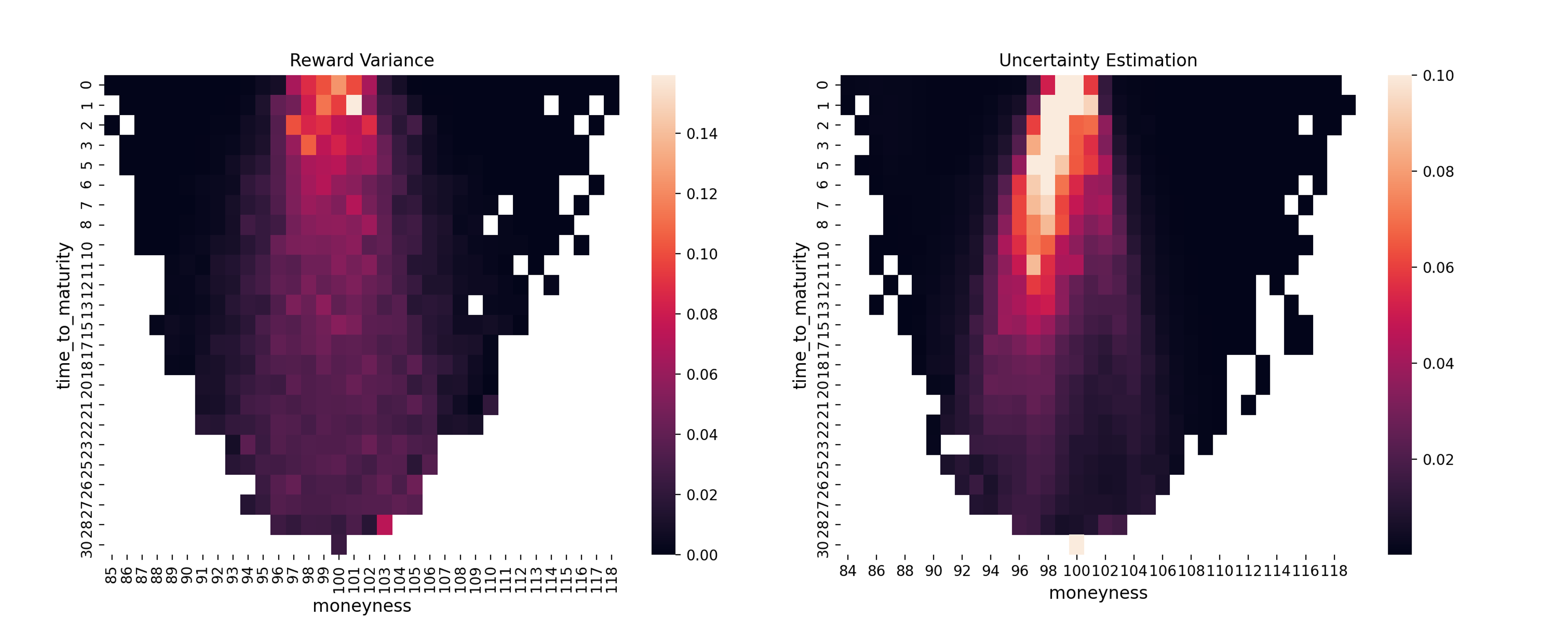}
	\end{center}
	\caption{\textbf{Uncertainty Estimation.} This figure shows the comparison of the uncertainty estimation of the model and the variance of the reward.}
	\label{uncertainty estimation}
\end{figure}

We consider the uncertainty estimation provided by the model and compare it with the real-world reward variance. Figure \ref{uncertainty estimation} shows the comparison of uncertainty estimation of the model and the variance of the reward. This is a heatmap that reflects the uncertainty of the returns. The brighter the part in the figure, the greater the uncertainty, while the darker the part, the smaller the uncertainty. The uncertainty distribution map reflects an obvious feature: when the stock price is near the strike price, the return uncertainty is the highest, and this feature becomes more pronounced as the expiration date approaches. This phenomenon is the same with the variance of reward, which is due to the greater volatility of option prices closer to the strike price.

\subsection{Real Data Application}
\subsubsection{Data Description And Preprocess}

\begin{figure}
	\begin{center}
	\includegraphics[width=0.95\columnwidth]{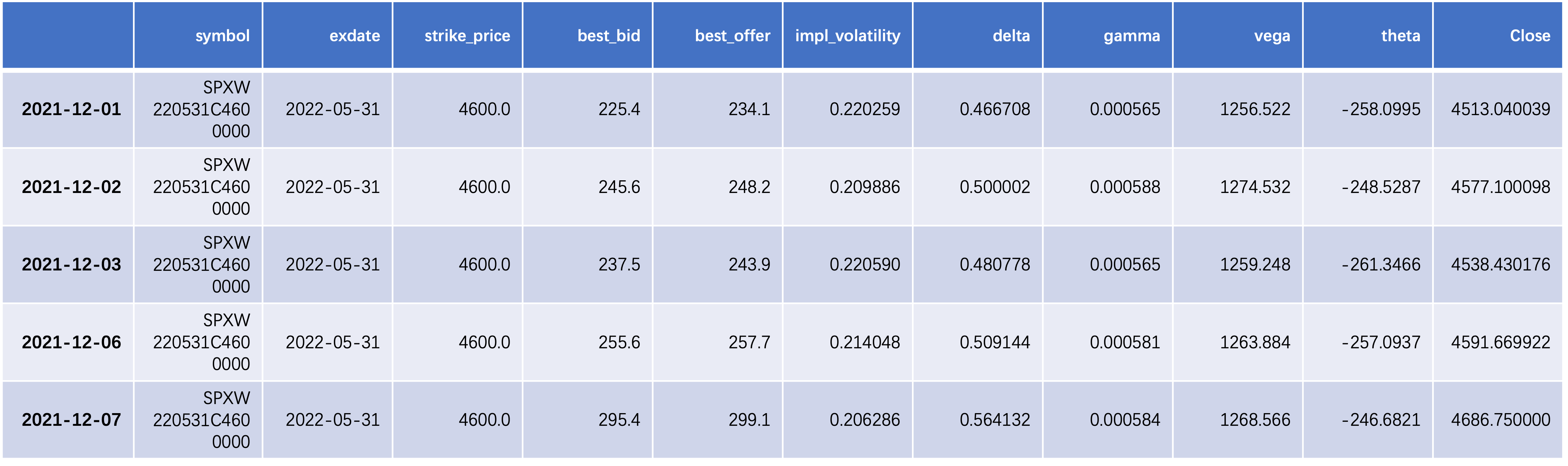}
	\end{center}
	\caption{\textbf{Option Data Sample.} This table provides a sample of option data, where the option symbol is SPXW 220531C4600000.}
	\label{Option Data Sample}
\end{figure}

\begin{figure}
	\begin{center}
	\includegraphics[width=0.95\columnwidth]{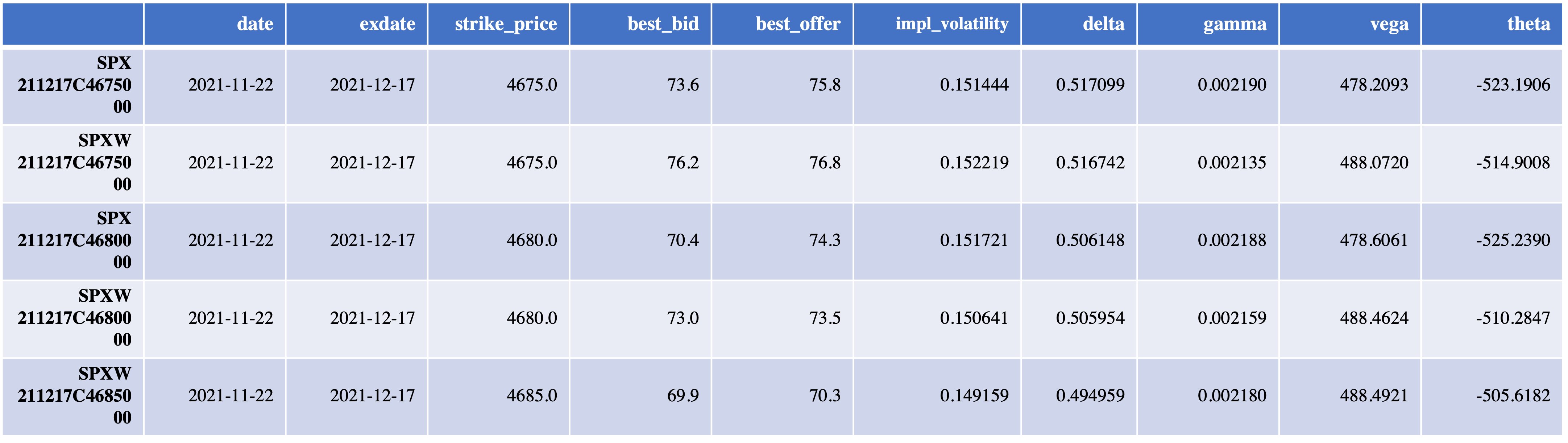}
	\end{center}
	\caption{\textbf{Option Chain Sample.} This table provides a sample of option chain with the different strike prices.}
	\label{Option Chain Sample}
\end{figure}

The data used in this study is options data for the S\&P 500, obtained from \cite{wachowicz2020wharton}. The data covers the period from January to December 2021 and includes 59,550 options with different expiration dates and strike prices. To focus on the optimal value of delta, the authors filtered out options with a delta value close to 0 or 1, i.e., options where the stock price deviates significantly from the strike price. Specifically, they focused on options with a stock price within $\pm$20\% of the strike price and an option length between 15 and 40 days. Additionally, the options data only provide the best bid and best ask prices, and not the actual option prices. To estimate the option prices, the authors used the mid-price of the best bid and best ask prices as an approximation. The option data sample can be seen in table \ref{Option Data Sample} and \ref{Option Chain Sample}.

For simulated data, they are concerned with the distribution of the total returns over the entire period, i.e., $\sum_t^T R_t$. However, for real-world options data, the lengths of different types of options are not the same, and therefore, the total sum of returns over the entire period may not be a good indicator of performance. To address this issue, the authors focus on the expected value and variance of the returns at each time step, rather than the total sum of returns over the entire period. This approach allows them to compare the performance of different options with different lengths and to analyze the risk and return characteristics of each option more accurately.

\begin{figure}
	\begin{center}
	\includegraphics[width=0.95\columnwidth]{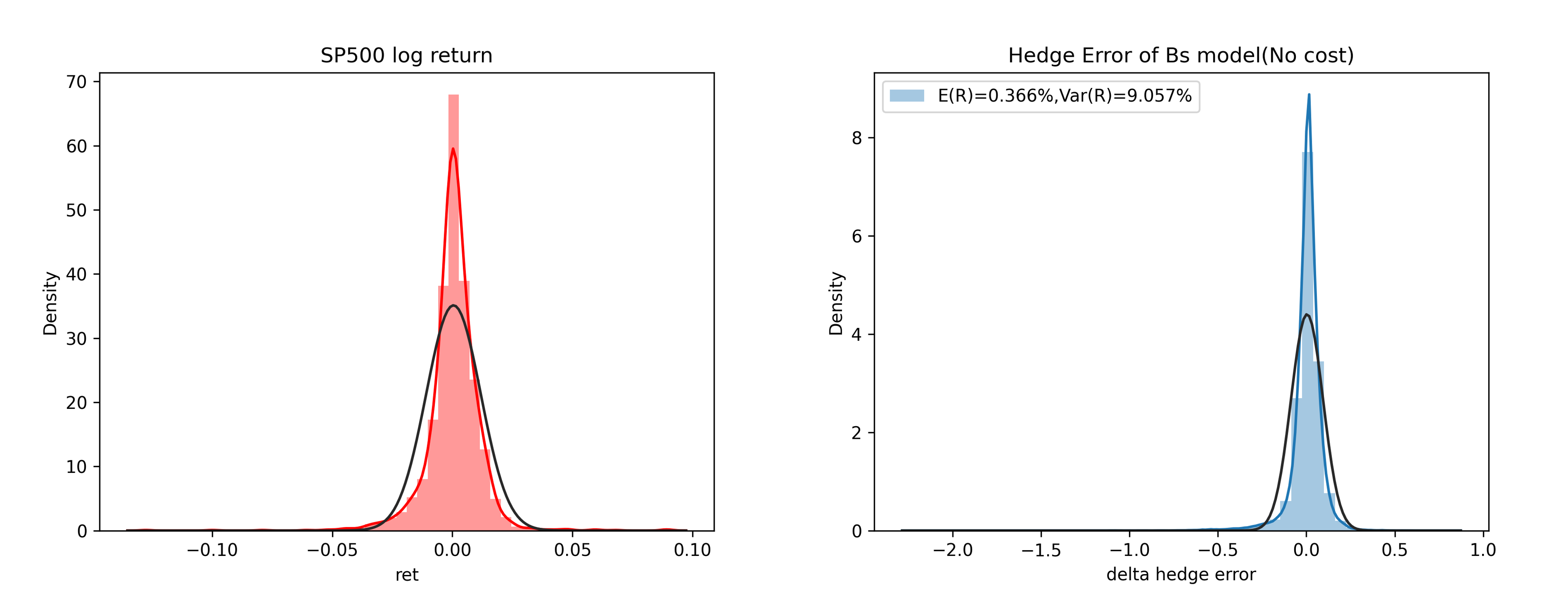}
	\end{center}
	\caption{\textbf{Fat tail challenge.} The fat-tail effect is a well-known phenomenon in financial data that poses a significant challenge for options traders, which violates the assumptions of the Black-Scholes model.}
	\label{Fat tail challenge}
\end{figure}

The fat-tail effect is a well-known phenomenon in financial data that poses a significant challenge for options traders. This phenomenon completely violates the assumptions of the Black-Scholes model, making it impossible for traders to rely solely on the theoretical delta formula. Figure \ref{Fat tail challenge} in the experimental analysis shows the distribution of log returns for the S\&P 500 index, as well as the distribution of hedging errors for the Black-Scholes model. Both graphs reflect the presence of the fat-tail effect, which is characterized by a higher probability of extreme events than would be predicted by a normal distribution.

\subsubsection{Feature Analysis}
\begin{figure}
	\begin{center}
	\includegraphics[width=0.8\columnwidth]{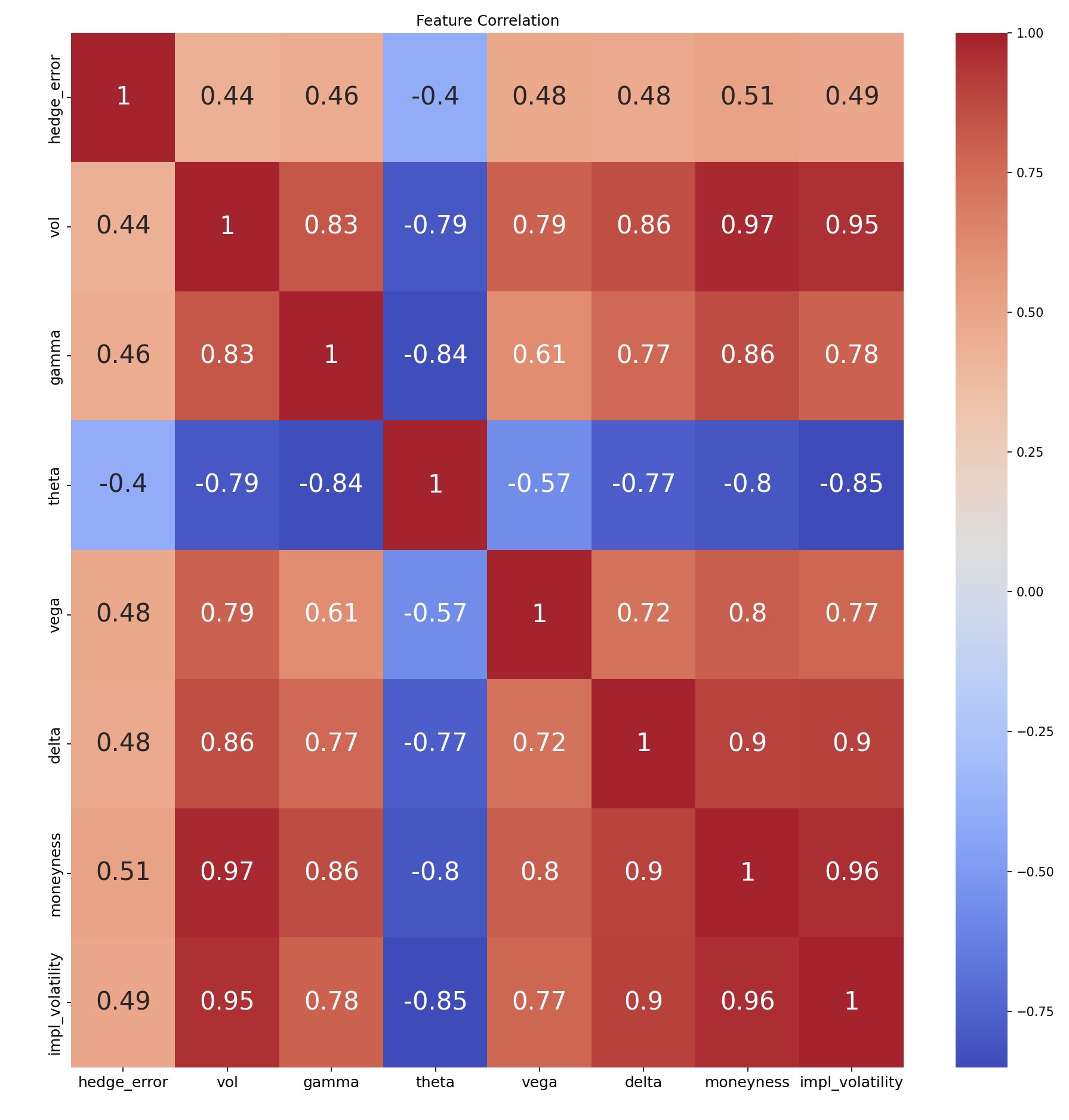}
	\end{center}
	\caption{\textbf{Feature Correlation.} This figure is a heat map of the relationship matrix between the residual term of the Black-Scholes model and multiple variables.}
	\label{Feature Correlation}
\end{figure}
To effectively solve the option hedging problem in real-world applications, it is crucial to extract relevant signals that can serve as input signals for the reinforcement learning model. Focusing on the cost problem of option dynamic hedging, the goal is to bring the value of $R_t$ closer to 0 by selecting an appropriate $N_{t+1}$. Given the baseline delta in the Black-Scholes model, we approach this as a boosting problem, as shown in \ref{residual}. By leveraging the residual learning technique, we aim to learn a residual function that can improve the accuracy of delta predictions and reduce the hedging error.

The term $C_t - C_{t+1} + \delta(S_{t+1}-S_t)$ in the formula is the residual term of the Black-Scholes model. The problem now becomes how to find $f(x_t)$ that can better fit the residual term. If we can find $f(x_t)$ such that formula \ref{residual} holds, then the delta of the dynamic hedging model is $\delta_{model} = \delta - f(x_t)$. If we view it as a regression problem, then $y = C_t - C_{t+1} + \delta(S_{t+1}-S_t)$, and the independent variable is $x = f(x_t) * (S_{t+1}-S_t)$.

\begin{equation}\label{residual}
    C_t-C_{t+1}+\delta(S_{t+1}-S_t) = f(x_t)(S_{t+1}-S_t).
\end{equation}

We provide a heat map of the relationship matrix between the residual term of the Black-Scholes model and multiple variables $x = f(x_t) * (S_{t+1}-S_t)$.

We choose the following features as state variables for the reinforcement learning model: time to maturity $\tau$, moneyness $S_t/K$, implied volatility of the stock $\sigma_t^{impl}$, Black-Scholes Vega, Black-Scholes Theta, Black-Scholes Gamma, 20-day historical volatility of the stock $\sigma_t^{20}$, and 30-day volatility of the stock $\sigma_t^{30}$.

\subsubsection{Results Analysis}
The empirical results of our analysis mainly compare the performance of the reinforcement learning model with that of the classical Black-Scholes model on real-world options data. We also focus on the improvement brought by uncertainty estimation in the training of the reinforcement learning model, as well as its potential application scenarios. Similar to the analysis of the simulation results, our evaluation of the real-world data focuses on the expected mean and variance of the option seller's profit. A higher expected value indicates that the hedging strategy is more effective and incurs less wear and tear, while a smaller variance indicates that the sale of options is less risky.

\begin{figure}
	\begin{center}
	\includegraphics[width=0.9\columnwidth]{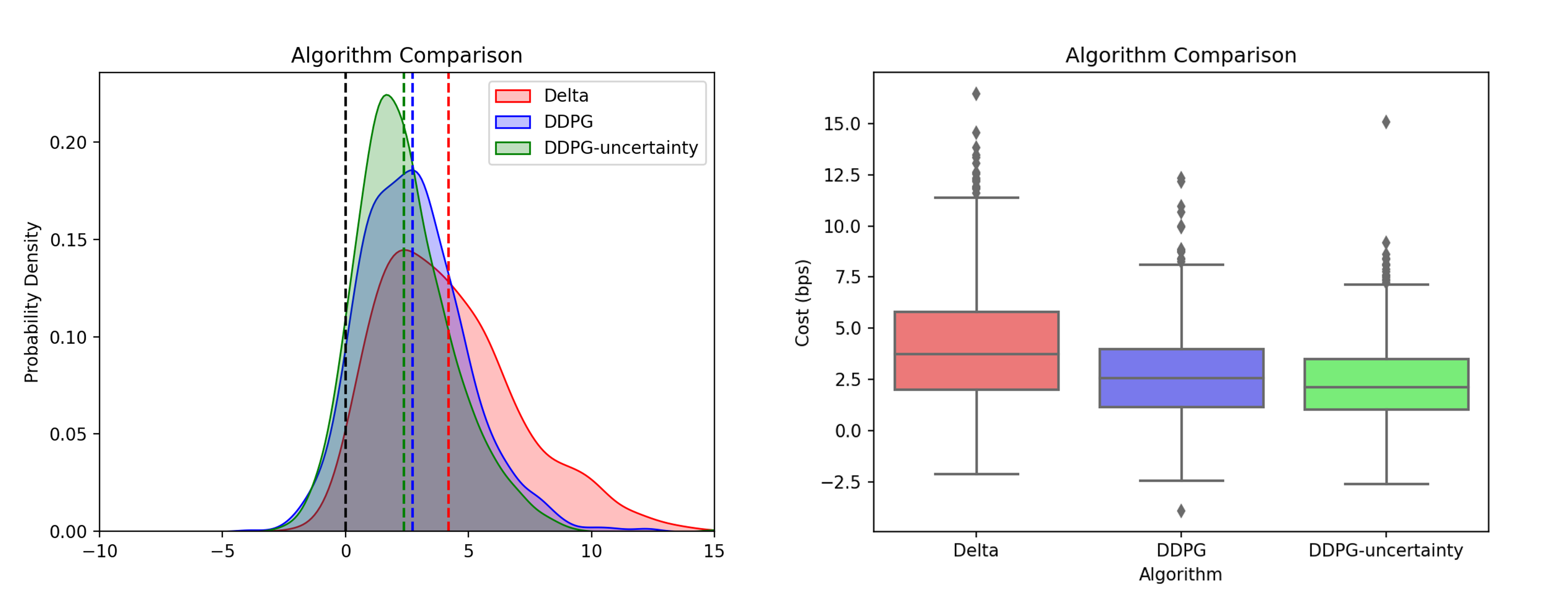}
	\end{center}
	\caption{\textbf{Strategy Comparison.} This figure compares the performance of different algorithms for dynamic hedging of S\&P 500 options. }
	\label{Strategy Comparison}
\end{figure}

The experimental analysis, presented in Figure \ref{Strategy Comparison} and Table \ref{tab:comparison_1}, compares the performance of different algorithms for dynamic hedging of S\&P 500 options. The experiment focuses on three different algorithms and demonstrates that the reinforcement learning model can more effectively manage wear and tear associated with hedging compared to the simple Black-Scholes model's $\delta$ value. This is evidenced by the significant improvement in the expected value metric. Furthermore, the results suggest that the reinforcement learning model can also reduce the risk associated with selling options. This finding was not evident in the analysis of the simulation data, which assumed constant volatility and was based on the basic assumptions of the Black-Scholes model. The real-world data used in the analysis is subject to greater volatility and uncertainty, which makes effective hedging more challenging. Nevertheless, the reinforcement learning model is able to adjust to changing market conditions and improve overall performance.

\begin{table}[ht]
\centering
\resizebox{0.8\textwidth}{!}{
\begin{tabular}{|c|c|c|c|}
\hline
 & $Mean(R_t)$ & $Var(R_t)$ & $Gain_{Delta}$ \\
\hline
Delta & -4.29 & 2.89 &  0\\
\hline
DDPG & -2.88 & 2.14 & 1.41 \\
\hline
$DDPG_{uncertainty}$ & -2.37 & 1.90 & 1.92 \\
\hline
\end{tabular}
}
\caption{\textbf{Algorithm Comparison.} This table provides a detailed comparison between three different strategies in terms of their mean and variance of the return.}
\label{tab:comparison_1}
\end{table}

The third algorithm involved in the experimental analysis is DDPG-uncertainty. In this approach, both aleatoric and epistemic uncertainties are simultaneously considered in the reinforcement learning training process. This allows for a more realistic and conservative trading strategy, given the uncertainty associated with stock price fluctuations. The experimental results presented in Table \ref{tab:comparison_1} show that DDPG-uncertainty can better control the wear and tear associated with dynamic hedging, as compared to algorithms that do not consider uncertainty. Additionally, DDPG-uncertainty is able to reduce risk to some extent, which is a key factor in effective risk management and improved returns.

In addition to examining the dynamic hedging effect, our analysis also focuses on the estimation of real-world uncertainty. To evaluate the level of uncertainty estimation, we use the variance of hedge error in real-world hedging scenarios as a measure of real-world uncertainty. Our goal is to provide accurate uncertainty estimates that can effectively estimate the uncertainty of the actual returns. Uncertainty in financial markets can arise from two main sources: first, the uncertainty of stock prices, as the rise and fall of stock prices and their magnitude are unpredictable (stock volatility); and second, the uncertainty of the model itself. Deep learning models are data-driven and rely on the assumption that the training and testing sets are independently and identically distributed (i.i.d.). However, this assumption may not always hold true in real-world applications. Therefore, providing uncertainty estimates along with model outputs is essential in the financial field. When the model provides a high level of uncertainty, traders can intervene in advance to avoid major trading accidents.

\begin{figure}
	\begin{center}
	\includegraphics[width=0.9\columnwidth]{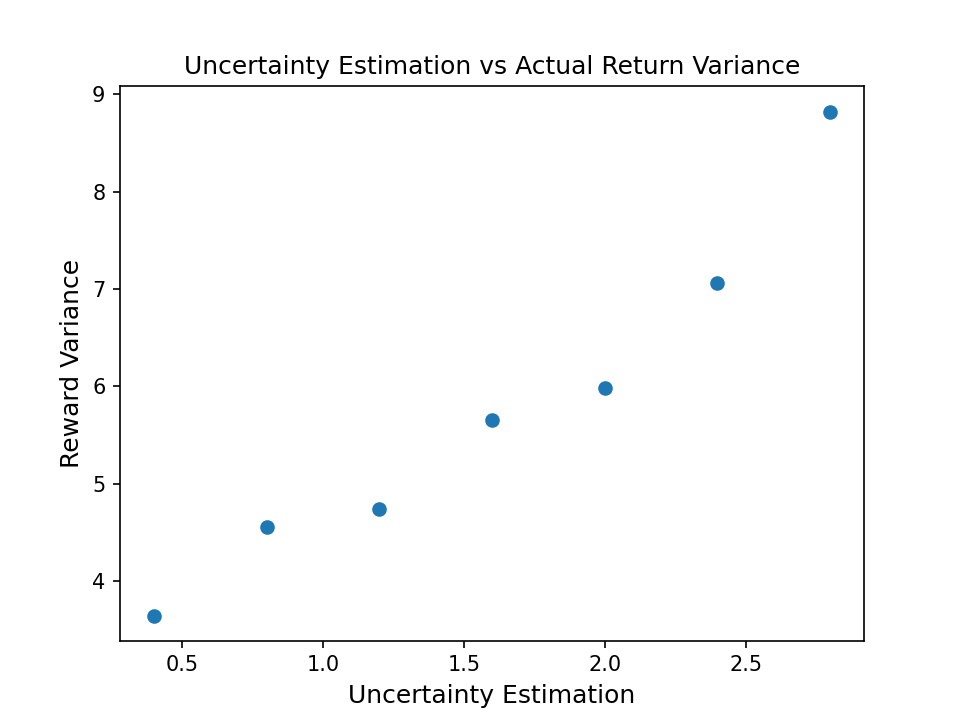}
	\end{center}
	\caption{\textbf{Correlation between Uncertainty Estimation and Actual Reward Variance.} This figure provides a visualization of the correlation between uncertainty estimates and the variance of corresponding real-world returns. We use the reward variance as an indicator to evaluate the accuracy of uncertainty estimation.}
	\label{Uncertainty reward}
\end{figure}

Figure \ref{Uncertainty reward} provides a visualization of the correlation between uncertainty estimates and the variance of corresponding real-world returns. Specifically, we randomly sampled 10,000 data points from the test set, each of which included two pieces of data: the uncertainty estimate provided by the model at the current step, and the actual reward obtained from hedging at the current step. We divided these samples into seven groups based on the value of the model's uncertainty estimate and calculated the variance of corresponding returns for each group. The results shown in the figure indicate a positive correlation between the uncertainty estimates provided by the model and the uncertainty of real-world returns. This suggests that the model's uncertainty estimates can effectively estimate the uncertainty of actual returns, and can provide traders with valuable information for risk management and decision-making.

\subsubsection{Examples}

\begin{figure}
	\begin{center}
	\includegraphics[width=0.8\columnwidth]{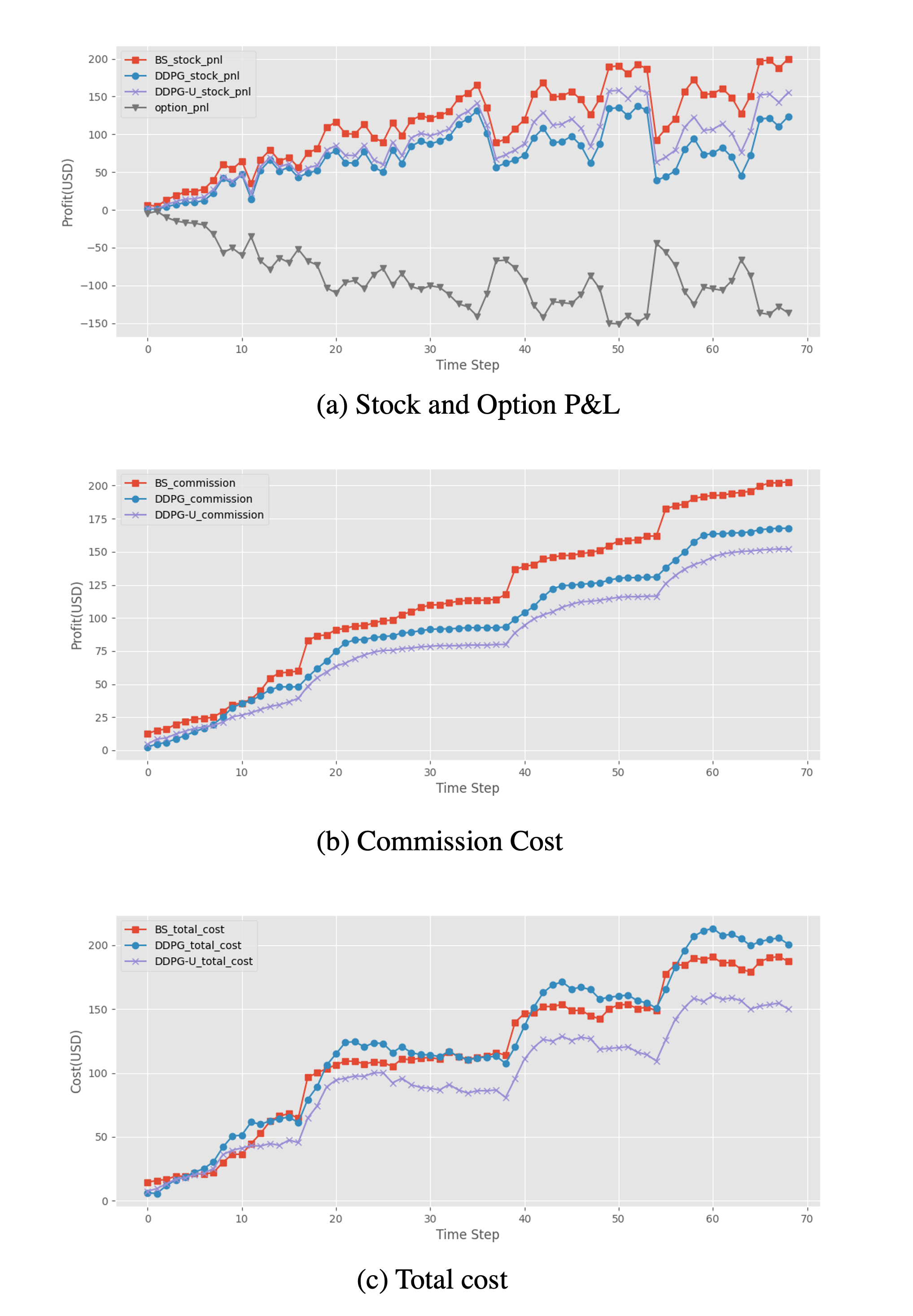}
	\end{center}
	\caption{\textbf{Strategy Comparison(Pnl decompose).} This figure provides a simple example of S\&P 500 option hedge process. \textbf{(a)}: option pnl represents the income generated by selling options, and stock pnl is the profit generated by buying and selling stocks during the hedging process. \textbf{(b)}: Commission is the cost generated by trading stocks. \textbf{(c)}: Cumulative Cost of the three hedging strategies.}
	\label{Strategy comparison(Pnl decompose)}
\end{figure}

\begin{figure}
	\begin{center}
	\includegraphics[width=0.9\columnwidth]{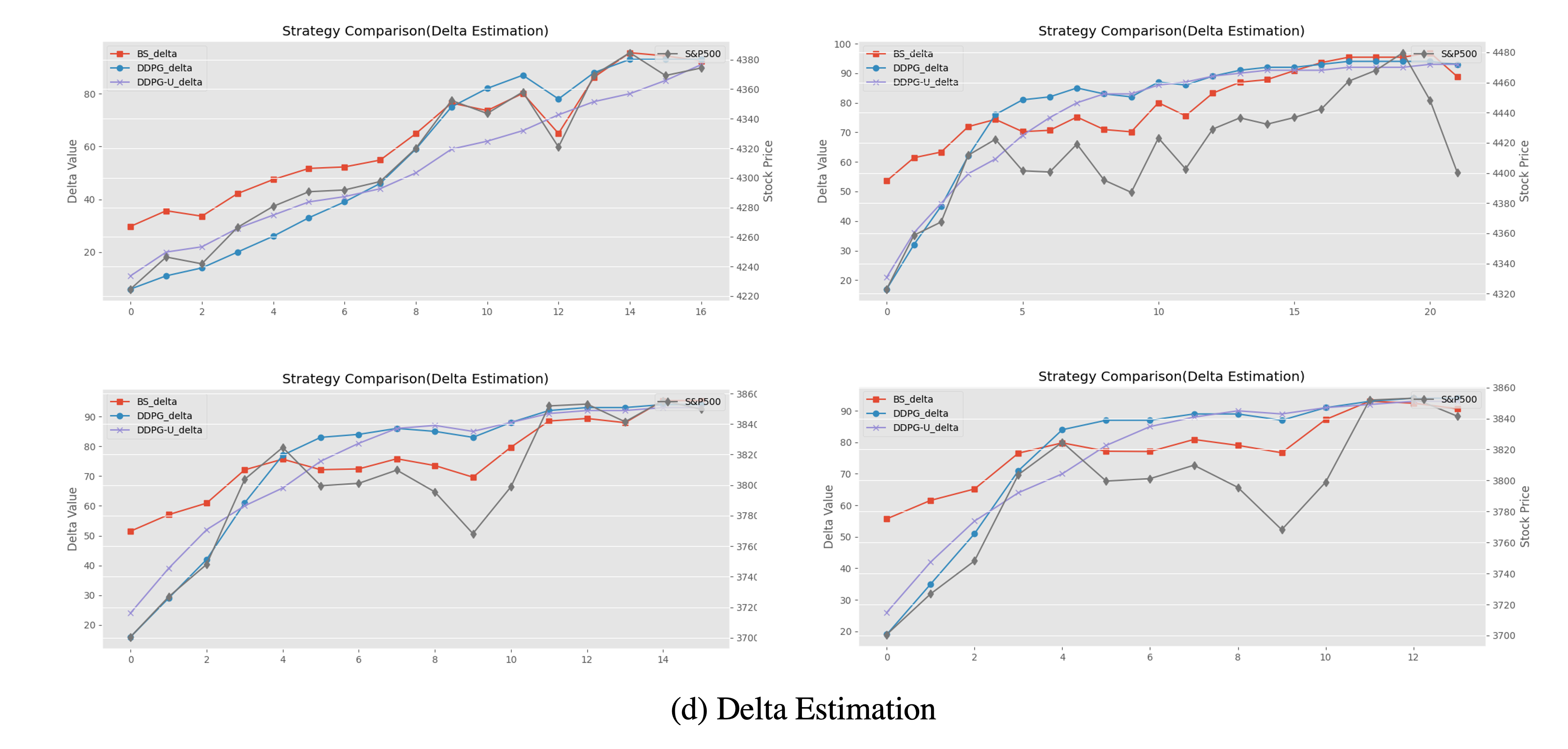}
	\end{center}
	\caption{\textbf{Strategy Comparison(Delta).} This figure provides four examples for the delta estimation of three different strategies.}
	\label{Strategy comparison(Delta)}
\end{figure}

To further analyze the decision-making process of the three dynamic hedging strategies, we provide a simple example of S\&P 500 options in Figure \ref{Strategy comparison(Pnl decompose)}. The figure shows the decomposition of the hedging process P\&L, with "BS" representing Delta hedging through the Black-Scholes model, "DDPG" representing the dynamic hedging model obtained through reinforcement learning, and "DDPG-U" representing the model obtained after considering uncertainty estimation during the model training process. 

The P\&L of the dynamic hedging process can be attributed to three main sources: income generated by selling options, profit generated by buying and selling stocks during the hedging process, and transaction costs incurred by buying and selling stocks. The dynamic hedging process involves buying and selling stocks to hedge against losses caused by option price fluctuations. As shown in the figure, all three strategies effectively hedge losses caused by selling call options through trading stocks (stock P\&L can offset the option P\&L). Additionally, trading stocks incur considerable wear and tear, with "DDPG-U" resulting in the least transaction cost wear and tear.

Figure c provides a comparison of the cumulative P\&L of the three hedging strategies, while Figure \ref{Strategy comparison(Delta)} specifically shows the delta decisions of the three strategies. By examining the cumulative P\&L chart, we observe that the P\&L changes of "DDPG-U" are relatively slow, and the cumulative P\&L value is also the largest, which is a desirable outcome. 

The delta estimation process provides further support for this result. As shown in Figure \ref{Strategy comparison(Delta)}, the delta value given by the Black-Scholes model is highly sensitive to fluctuations in stock prices, especially at the timestamp when the option is just sold. This can cause the stock position to quickly reach a relatively high level, leading to significant wear and tear from buying and selling stocks. In contrast, "DDPG" provides smoother up and down fluctuations in the Delta value, although it is still slightly larger than "DDPG-U". "DDPG-U" relies on relatively smooth position adjustments, making the entire dynamic hedging process more natural and effective.

\section{Conclusion}

This work demonstrates the application of reinforcement learning methods to solve the problem of dynamic hedging of options, providing greater flexibility compared to traditional analytical solutions. Real-world financial applications often involve more complex problems than the modeling assumptions of theoretical models. With the advent of deep learning, and specifically the development of reinforcement learning models, it is possible to model sequential decision-making problems using deep models. To address the uncertainty issues commonly present in financial problems, we have incorporated the effects of both aleatoric and epistemic uncertainties in our modeling.

Our contribution to the field of dynamic hedging of options mainly focuses on the following aspects. Firstly, we have developed a reinforcement learning-based interactive environment, in which the agent optimizes the action space for dynamic hedging of options using rewards provided by the environment. This approach provides greater flexibility compared to traditional analytical solutions and enables traders to adapt to complex and dynamic market conditions. Secondly, we have made certain adjustments to the Q function to measure the risk of dynamic hedging more effectively. Specifically, we have incorporated the effects of transaction costs and wear and tear in the modeling process, making it more convenient for risk measurement and enabling traders to better manage risk exposure. Finally, to deal with the unpredictable fluctuations in financial markets, we have fully modeled aleatoric and epistemic uncertainties in the modeling process. This ensures that the trading strategy is not overly optimistic and can effectively deal with unexpected situations in real-world scenarios.

\bibliographystyle{unsrtnat}
\bibliography{references}  






\end{document}